\documentclass[a4paper,10pt,reqno]{amsart}
\usepackage{verbatim}
\usepackage{amssymb}
       
\usepackage{enumerate}   
\usepackage[active]{srcltx}   
\numberwithin{equation}{section}        
\usepackage{graphicx}
\usepackage[inline]{showlabels}
       
\usepackage{t1enc}       
\usepackage[utf8x]{inputenc}        
       
\usepackage{tikz}

\newtheorem{theorem}{Theorem}[section]        
\newtheorem{lemma}[theorem]{Lemma}       
       
\newtheorem{problem}[theorem]{Problem}       
       
\newtheorem{proposition}[theorem]{Proposition}

\newtheorem{claim}{Claim}[theorem]

\theoremstyle{definition}       
\newtheorem{definition}[theorem]{Definition}       
\theoremstyle{remark}       
         
\newcommand{\mc}[1]{\mathcal{#1}}       
\newcommand{\mbb}[1]{\mathbb{#1}}

\newcommand{\setm}{\setminus}       
\newcommand{\empt}{\emptyset}       
\newcommand{\subs}{\subset}

\def\<{\left\langle}       
\def\>{\right\rangle}       

\newcommand{\ddeg}{\operatorname{d}}

\newcommand{\dgl}{\ddeg^{<}}

\newcommand{\ccc}{\operatorname{c}}
\newcommand{\CCC}{\operatorname{C}}

\newcommand{\CC}[2]{\CCC_{#2}(#1)}

\newcommand{\cc}[2]{\ccc_{#2}(#1)}

\newcommand{\www}{\operatorname{w}}
\newcommand{\WWW}{\operatorname{W}}

\newcommand{\WW}[2]{\WWW(#2,#1)}

\newcommand{\ww}[2]{\www(#1,#2)}

\newcommand{\gmin}[1]{\min #1}
\newcommand{\gmax}[1]{\max #1}

\newcommand{\isf}{\mc F_{{[{\omega}]}^{2}}}

\author[T. Csern\'ak]{Tam\'as Csern\'ak}
\address
      {Eötvös University of Budapest, Hungary  }
\email{tamas@csernak.com}

\author[L. Soukup]{Lajos Soukup}
\address
      { Alfr{\'e}d R{\'e}nyi Institute of Mathematics, Eötvös Loránd Research Network, Budapest, Hungary  }
\email{soukup@renyi.hu}

\title%[Elusive properties]%
   {Elusive properties of infinite graphs}
\thanks{The preparation of this paper was partially
supported by  OTKA grant K 129211}

\date{\today}

\begin{document}

\begin{abstract} 
      A graph property  
      is said to be {\em  elusive ({\rm or} evasive)} if every algorithm testing this property by asking questions of the 
      form "{\em is there an edge between vertices $x$ and $y$}'' requires, in the worst case, 
      to ask about all pairs of vertices.

      The unsettled Aanderaa–Karp–Rosenberg conjecture is that 
      every non-trivial monotone graph property is elusive for finite vertex sets.

We show that the situation is completely different for infinite  vertex sets: 
the monotone graph properties ``{\em every vertex has degree at least $n$}'' and 
      ``{\em every connected component  has size at least $m$}'', where $n\ge 1$ and 
      $m\ge 2 $ are  natural numbers, 
      are not elusive
      for infinite vertex sets, but        
      the monotone graph property ``{\em the graph contains a cycle }'' is elusive
      for arbitrary vertex set.
      
  On the other hand, we also prove that every algorithm testing  some natural 
  monotone graph properties, e.g  
  ``{\em every vertex has degree at least $n$}''  or  ``{\em connected}''  on the vertex set ${\omega}$    
  should check "lots of edges", more precisely, 
  all the edges of an infinite complete subgraph.
    \end{abstract}       
        
\maketitle

\section{Introduction}  

Given a graph property $R$ and a vertex set $V$,  let us consider the following 
game between two players, Alice, the seeker, and Bob, the hider.
 First the hider takes a graph
$G$ with vertex set $V$. Then, in each move of the game, the seeker asks
the hider whether a certain edge $e$ is in $G$ or not. The game terminates
when the seeker can decide whether $G$ has property $R$.

Rosenberg \cite{Ro} conjectured that for every non-trivial graph property $R$
and for each finite vertex set $V$ the seeker should ask $\Omega(|V|^2)$ many edges in the worst
case. In \cite{LBB} this conjecture was disproved by showing that the 
"{\em being a scorpion graph}"
property can be decided by the seeker in $O(|V|)$ steps.

However, a weaker conjecture of Aanderaa and Rosenberg was proved by 
Rivest and Vuillemin in \cite{RiVu}: if $R$ is a non-trivial {\em monotone} 
graph property, then  for each finite vertex set $V$ the seeker should ask $\Omega(|V|^2)$ 
many edges in the worst case; in this context, 
a property is monotone if it remains true when edges are added.

A graph property $R$ is said to be {\em elusive on vertex set $V$}
if the hider has a strategy such that the seeker needs query all the possible edges to decide 
whether $G$ has property $R$ or not.

A stronger, unsettled version of  the Aanderaa–Rosenberg conjecture, called the 
evasiveness conjecture or the Aanderaa–Karp–Rosenberg conjecture, states 
that  
every non-trivial monotone graph property is elusive on each finite vertex set $V$.

Many partial results were proved in connection with this conjecture,  
see e.g. \cite{Bo}, \cite{LBB}.  
Kahn, Saks and Sturtevant  \cite{KaSaSt84} proved it in the case when $|V|$ is a prime power, and 
Yao \cite{Ya88} proved an analogous conjecture for bipartite graphs.  

In this paper we %investigate these problems for infinite vertex sets to 
show that the situation is quite different when the vertex set is infinite.
 
First, in section  
\ref{sc:elusive} we prove that the following natural monotone graph properties are not elusive 
for infinite vertex sets: 
\begin{enumerate}[(i)]
\item[$(D_n)$] "{\em the degree of each vertex is at least $n$}", 
\item[$(C_m)$] "{\em the connected components have  size at least $m$}",  
\end{enumerate}
for each natural number $n\ge 1$ and $m\ge 2$
(see Theorem \ref{tm:gen}).

On the other hand, we also prove that the monotone graph property 
"{\em the graph contains a cycle}" is elusive for arbitrary vertex sets 
(see Theorem \ref{tm:cycle}).

At the moment we do not have  any reasonable conjecture to characterize 
(large classes of) monotone graph properties which 
are elusive for infinite vertex sets.

In section \ref{sc:weak-elusive} we formulate a theorem for infinite vertex sets 
which 
has the flavor of
the result of Rivest and Vuillemin  from \cite{RiVu}: 
although  the properties $D_n$ and $C_m$ above are not elusive, but 
Alice, the seeker should check "lots of edges" 
before she can decide if a given  graph has that property (see Theorem \ref{tm:we-degree}). 

We also show that the "{\em being a scorpion graph}"
property can be decided by the seeker asking "few" edges (see Theorem \ref{tm:scorpion}).

\section{Basic notions  and definitions}

First, we give  a game-theoretic reformulation of the problem 
which also  covers e.g. the bipartite case.

Based on  the terminology of Bollobás \cite{Bo},
given a graph $H=\<V,E^*\>$ we say that  a triple $G = \<V, E, N\>$ is an 
 {\em $H$-pregraphs} iff 
 $E \cup  N\subs E^*$ and 
$E \cap N = \empt$. We will say  that $V$ is the set of {\em vertices}, $E$ is the set of {\em edges}, and $N$ is
the set of {\em nonedges}, $P = E \cup  N$ is the set of {\em determined pairs}    of $G$,
and $U = E^*\setm P$ is the set of {\em undetermined pairs}    of $G$.
We will write $G=\<V_G,E_G,N_G\>$, $P_G=P$, $U_G=U$, 
 $\gmin G = \<V_G,E_G\>$ and $\gmax G=\<V_G, E_G\cup U_G\>$. 

 We say that a graph $G=\<V,E\>$ is an \emph{$H$-graph} iff 
 $E\subs E^*$. 
 If $G=\<V,E\>$ is a graph, 
the set of {\em nonedges} of $G$ is 
$N_G={[V]}^{2}\setm E$.

If $H$ is the complete graph on $V$, then we will write ``{\em pregraph on $V$ }''
or simply  ``{\em pregraph}'' instead of ``$H$-pregraph''. 

If $G_1$ is a pregraph on $V$, and $G_2$ is either a graph or a pregraph on $V$, 
 we say that $G_2$ 
{\em extends}  $G_1$, $G_1\le G_2$, 
provided  $E_{G_2}\supseteq E_{G_1}$ and $N_{G_2}\supseteq N_{G_1}$.

%We will  write \emph{extension } instead of $H$-extension, if $H$ is clear from the context.  

Given a graph property $R$ we say that an $H$-pregraph $G=\<V,E,N\>$ has property $R$
iff the graph $\gmin G=\<V,E\>$ has property $R$.

\medskip

Since we want to consider the elusiveness of properties of bipartite graphs as well, we should restrict the set of possible edges. 
That is the reason that we need an addition parameter $E^*$ in the next definition: 
the set of allowed edges.  

\begin{definition}\label{df:w-elusive1} 

      Let $R$ be a graph property, 
      $V$ be the set of vertices, and $E^*\subs {[V]}^{2}$
      be the set of 
      {\em  allowed edges}. 
      Write $H=\<V,E^*\>$.

      Define the game $\mbb E_{H,R}$ between two players, Alice and Bob, as follows: 
      \begin{enumerate}[(1)]
      \item During the game Alice and Bob construct a (transfinite) sequence $\<G_{\alpha}:{\alpha}\le {\beta}\>$ of 
      $H$-pregraphs with  $G_0=\<V,\empt,\empt\>$.
      \item  
      The game terminates when either every $H$-graph extenting  $G_{\beta}$ 
      has property $R$, or no $H$-graph extenting $G_{\beta}$ has property $R$.
      \item If ${\alpha}$ is a limit ordinal, then let 
      $E_{G_\alpha}=\bigcup_{{\zeta}<{\alpha}}E_{G_\zeta}$ and 
      $N_{G_\alpha}=\bigcup_{{\zeta}<{\alpha}}N_{G_\alpha}$.
      \item If ${\alpha}={\gamma}+1$ and the game has not terminated, then  $G_{\gamma}$ has undetermined pairs, and  
      \begin{enumerate}[(i)]
      \item Alice picks an undetermined pair $e_{\gamma}\in U_{G_{\gamma}}$,
      \item Bob decides if  $e_{\gamma}$ is an edge, or a nonedge in $G_{\alpha}$, i.e. Bob 
      selects an $H$-pregraph $G_{\alpha}$ such that $G_{\alpha}>G_{\gamma}$ and 
      $P_{G_\alpha}=P_{G_\gamma}\cup \{e_{\gamma}\}$. 
      \end{enumerate}
      \item Bob {\em wins} iff 
      $P_{\beta}=E^*$, i.e. $\{e_{\alpha}:{\alpha}<{\beta}\}=E^*$ and so $U_{\beta}=\empt$.
      \end{enumerate}
      If the property $R$ is monotone, then (2) can be written in the following form:
      \begin{enumerate}[(1')]\addtocounter{enumi}{1}
      \item 
      The game terminates in turn ${\beta}$ if either  
      $\gmin{G_{\beta}}=\<V,E_{\beta}\>$ has property $R$ or 
      $\gmax{G_{\beta}}=\<V,E_{\beta}\cup U_{\beta}\>$
      does not have property $R$.
      \end{enumerate}   
To simplify the notation, we will write $E_{\alpha}$ for $E_{G_{\alpha}}$, 
$N_{\alpha}$ for $N_{G_{\alpha}}$, $P_{\alpha}$ for $P_{G_{\alpha}}$,  
and $U_{\alpha}$ for $U_{G_{\alpha}}$.   

      If we do not have any restriction on edges, i.e. if
$H=K_V=\bigl\langle V,{[V]}^{2}\bigr\rangle$, 
 then we write $\mbb E_{V,R}$ for $\mbb E_{K_V,R}$. 

      \end{definition}

\begin{definition} 
      A graph property $R$ is {\em elusive in a graph $H$} iff Alice does not 
      have a winning strategy in the game $\mbb E_{H,R}$.
      A graph property $R$ is {\em elusive on a set $V$} iff Alice does not 
have a winning strategy in the game $\mbb E_{V,R}$.

A graph property $R$ is {\em strongly elusive in a graph $H$} iff Bob has a winning strategy in the game 
$\mbb E_{H,R}$.
A graph property $R$ is {\em strongly elusive on a set $V$} iff Bob has a winning strategy in the game $\mbb E_{V,R}$.

\end{definition}

\begin{definition}
For natural numbers $n\ge 1$ and $m\ge 2$, let $D_n$ and $C_m$ denote the  
following graph properties: 
\begin{enumerate}[(i)]
      \item[$(D_n)$] "{\em the degree of each vertex is at least $n$}", 
      \item[$(C_m)$] "{\em the connected components have  size at least $m$}",  
      \end{enumerate}
\end{definition}

Our notation is standard, see e.g. \cite{Di-Book}. Given sets  $X$ and $Y$, let 
$[X,Y]=\big\{\{a,b\}:a\in X,\ b\in Y, a\ne b\big\}$.  If $G=\<V,E\>$ is a graph and $v$ is a vertex, 
 $\ddeg_G(v)$ denotes 
the degree of the vertex $v$, $\CC{v}{G}$ denotes the connected component of $v$, and 
$\cc{v}{G}=|\CC{v}{G}|$. 
 If $V$ is clear from the context, we will write $\ddeg_E(v)$
for $\ddeg_G(v)$, $\CC{v}{E}$ for $\CC{v}{G}$ and $\cc{v}{E}$ for $\cc{v}{G}$.

\section{Elusive and non-elusive properties}\label{sc:elusive}

First we show that  
some natural monotone graph properties are elusive. 
\begin{theorem}\label{tm:cycle}
      The  monotone graph property   $R$ 
      \begin{displaymath}
      \text{``{\em the graph contains a cycle }''}
      \end{displaymath}
       is strongly 
      elusive
      on any  vertex set $V$ with $|V|\ge 3$. 
      \end{theorem}

      Instead of Theorem \ref{tm:cycle} we prove the following stronger result.

\begin{theorem}\label{tm:cycle-gen}
Let $H=\<V,E^*\>$ be a graph and  let $R$  denote the  monotone graph property   
\begin{displaymath}
\text{``{\em the graph contains a cycle }''}
\end{displaymath}    
Then following  statements are equivalent:
\begin{enumerate}[(1)]
\item every connected component of $H$ is $2$-edge connected,
\item $R$
is elusive in $H$, i.e. Alice  does not 
have a winning strategy in the game $\mbb E_{H,R}$
\item $R$
 is strongly 
elusive in $H$, i.e.  Bob has a winning strategy in the game 
$\mbb E_{H,R}$.
\end{enumerate}
      \end{theorem}

      \begin{proof}[Proof of Theorem \ref{tm:cycle-gen}]
 (2) implies (1).

We show that $\neg(1)$ implies $\neg(2)$.  Assume that $K$ is a connected component of $H$, but $K$ is not $2$-edges connected,
i.e. there is an edge $e\in E^*$ such that $K\setm \{e\}$ is not connected.

Then Alice has the following winning strategy in the game $\mbb E_{H,R}$:
Alice enumerates all the edges of $H$ as $\{e_{\alpha}:{\alpha}\le {\beta}\}$ with $e_{\beta}=e$, and she asks 
$e_{\alpha}$ in the ${\alpha}$th step. Then $\gmin {G_{\beta}}=\<V,E_{\beta}\>$ contains a cycle iff 
$\gmax{G_{\beta}}=\<V,E_{\beta}\cup\{e\}\>$ contains a cycle because  
there is no 
cycle in $H$ which contains $e$. So  the game terminates after at most ${\beta}$ steps 
and $U_{\beta}=\{e\}\ne \empt$. Thus, Alice wins. 

\medskip\noindent
 (3) implies (2).  Trivial.

\medskip\noindent
 (1) implies (3).

 We can assume that $H$ is connected because Bob can image that he plays a separate game
 in each connected component of $H$, and if he wins in each component, then he also 
 wins in $\mbb E_{H,R}$.
 
 So from now on we assume that $H$ is connected. 
      Let $V=V(H)$. 
      We show that the following greedy  algorithm gives a winning strategy for Bob in 
      $\mbb E_{H,R}$.
      
      Assume that in the ${\alpha}$th turn of the game we have a 
      pregraph $G_{\alpha}=\<V_{\alpha},E_{\alpha},N_{\alpha}\>$ and 
      Alice picked the pair $e_{\alpha}\in U_{\alpha}=E^*\setm (E_{\alpha}\cup N_{\alpha})$. 
      
      If $\<V,E_{\alpha}\cup \{e_{\alpha}\}\>$ is cycle-free, then Bob declares that $e_{\alpha}$ is en edge,
      i.e. $E_{{\alpha}+1}=E_{\alpha}\cup \{e_{\alpha}\}$. Otherwise, $e_{\alpha}$ will be a nonedge.
      
      Assume that the game terminates after ${\beta}$ turns.
      \begin{lemma}\label{lm:cfree}
      For each ${\alpha}\le {\beta}$, the graph $\<V,E_{\alpha}\>$ is cycle-free.
      \end{lemma}
      
      \begin{proof}
      Trivial by transfinite induction.
      \end{proof} 
      
      \begin{lemma}\label{lm:stable}
            For each ${\alpha}\le {\beta}$, if $\<V,E_{\alpha}\>$ is connected and $U_{\alpha}\ne \empt$, then 
            $\gmax{G_{\alpha}}=\<V,E_{\alpha}\cup U_{\alpha}\>$ contains a cycle. 
      \end{lemma}

      \begin{proof}
       Indeed, if  $\<V,E_{\alpha}\>$ is connected
        and $e\in {[V]}^{2}\setm E_{\alpha}$ then $\<V,E_{\alpha}\cup\{e\}\>$ contains a cycle.
      \end{proof}
      
      \begin{lemma}\label{lm:triangle}
            For each ${\alpha}\le {\beta}$,      if $\<V,E_{\alpha}\>$ is not connected then 
            $\gmax {G_{\alpha}}=\<V,E_{\alpha}\cup U_{\alpha}\>$ contains a cycle.
      \end{lemma}
      
      \begin{proof}
            Let  $\<A,B\>$ be a partition of $V$ such that $[A,B]\cap E_{\alpha}=\empt$.
      Then there is no nonedge between $A$ and $B$.
            
      Pick $a\in A$ and $b\in B$. Since $H$ is 2-edge-connected, there are two 
      edge disjoint paths between $a$ and $b$, $e_0\dots e_n$ and
      $f_0\dots f_m$ in $H$. 

If $e_i$ is an edge inside a connected component of $\<V,E_{\alpha}\>$, then replace 
$e_i$ with a path in $E_{\alpha}$ between the endpoints of $e_i$.

Since there is no non-edge between distinct connected components of $\<V,E_{\alpha}\>$,
 we obtain the path $\bar e'=e'_0\dots e'_k$ between $a$ and $b$ in 
$\gmax{G_{\alpha}}=\<V,E_{\alpha}\cup U_{\alpha}\>$.

Similarly, we obtain the path $\bar f'=f'_0\dots f'_\ell$ between $a$ and $b$ in $\gmax{G_{\alpha}}$
from $f_0\dots  f_m$. 

Then  $\bar e'$ contains at least one edge $e'$ between $A$ and $B.$
Since $e'$ is not in $\bar f'$, there is a path 
$g_0\dots g_s$ between the endpoints of $e'$ containing edges from  $\bar e\cup \bar f\setm \{e'\}$.

Thus, $e'g_0\dots g_s$ is a cycle in $\gmax{G_{\alpha}}=\<V,E_{\alpha}\cup U_{\alpha}\>$.     
      \end{proof}
      
      By lemma \ref{lm:cfree}  $\gmin G_{\beta}=\<V,E_{\beta}\>$ does not contain a cycle.
      Since the game has terminated, $\gmax{G_{\beta}}=\<V, E_{\beta}\cup U_{\beta}\>$ does not contain a cycle either.   
      Thus, by lemma  \ref{lm:triangle} the graph $G_{\beta}$ is connected and so  
      by Lemma \ref{lm:stable}  
      $U_{G_{\beta}}=\empt.$ 
      
      So Bob wins. 
      \end{proof}

\medskip

The monotone graph property "{\em every vertex has infinite degree}" is clearly not elusive 
on infinite vertex sets. 
However, next  we show that even certain natural, not "tailor made", monotone graphs properties 
are also not elusive.
First we should recall some definition.

\begin{definition}
      For a cardinal ${\kappa}$ and a natural number $n\ge 2$, the  {\em Turán graph 
      $T_{{\kappa},n}$}  is defined  
      as follows: $V(T_{{\kappa},n})={\kappa}\times n$
      and $E(T_{{\kappa},n})=\big\{\{\<{\alpha},i\>,\<{\beta},j\>\}:{\alpha},{\beta}\in {\kappa}, i\ne j<n\big\}$. 
      \end{definition}
      
Let us remark that $T_{{\kappa},2}$ is just the complete balanced bipartite graph of cardinality ${\kappa}$. 

\begin{definition} The 
      {\em Cantor graph} $C$ is defined as follows: its vertex set is the set of all finite
      0-1 sequences, and $\{s,t\}$  is an edge iff $s\subsetneq t$ or $t\subsetneq s$. 
      \end{definition}

\begin{theorem}\label{tm:gen}
If 
\begin{enumerate}[(i)]
\item $P$ denotes either the graph property $D_n$ for some $n\ge 1$ or
the graph property $C_m$ for some $m\ge 2$, and
\item $H$ is either  
the complete graph $K_{\kappa}$ 
or the Turan graph $T_{{\kappa},k}$ for some infinite cardinal ${\kappa}$ and $k\ge 2$,
or  the Cantor graph,
\end{enumerate}
 then  Alice has a winning strategy in the game $\mbb E_{H,P}$, i.e. 
\begin{center}
      the property $P$ is not elusive in $H$.
\end{center}
\end{theorem}

The previous theorem actually consists of  6 statements, but as we will see,  one 
can find a common generalization of all of them   (see Theorem \ref{tm:braided-w} below). 

First we find some  property shared by  the complete graphs, the Turan graphs and 
the Cantor graph (see  Definition \ref{df:braided} below), then 
we isolate some properties of $C_m$ and $D_n$ which makes possible to prove
the failure of elusiveness.

\begin{definition}\label{df:braided}
Let $H=\<V,E^*\>$ be a graph. A vertex set $L\subs V$ is a {\em covering} set iff
for each $v\in V\setm L $ there is $a\in L$ with $\{v,a\}\in E^*$.

We say that $H$ is {\em braided} iff it contains ``lots of'' finite covering sets: 
for each $W\in {[V]}^{<|V|}$ there is a 
finite covering set $L\in {[V\setm W]}^{<{\omega}}$.
\end{definition}

\begin{proposition}\label{pr:braided}     
 Given an infinite cardinal ${\kappa}$,      
the infinite complete graph $K_{\kappa}$, %the balanced bipartite graph $K_{{\kappa},{\kappa}}$,
the Turan graphs $T_{{\kappa},k}$ for $2\le k<{\omega}$, and the Cantor graph $C$ are braided.  
\end{proposition}

\begin{proof}
The finite set  of 0-1 sequences of length  $k$ is a covering set in $C$ for $k\in {\omega}$. 
So $C$ is braided. The other statements from this proposition are trivial.   
\end{proof}

\begin{definition}\label{df:w}
      Assume that $V$ is an infinite set, $H=\<V,E^*\>$ is a graph,  $n\in {\omega}$, and   
       $w:V\times \mc P(E^*)\to n+1$ is a function.
      We say that $w$ is
\begin{enumerate}[(w1)]
 \item \label{non-trivial} {\em non-trivial} iff $w(a,\empt)=0$ and $w(a,E^*)=n$
 for each $a\in V$;     
\item  \label{w:deg} {\em degree restricted} iff there is $K\in{\omega} $
such that $\ddeg_E(a)\ge K$ implies $w(a,E)=n$ for each $E\subs E^*$ and $a\in V$;
\item \label{w:finitely determined}   {\em finitely determined} 
iff $w(a,E)=\max\{w(a,E'):E'\in {[E]}^{<{\omega}}\}$ for each $E\subs E^*$ and $a\in V$;
\item \label{w:bounded} {\em bounded} iff there is $M\in {\omega}$
and there is a function $W:V\times \mc P(E^*)\to {[V]}^{\le M}$ such that
\begin{enumerate}[(i)]
\item  $a\in W(a,E)$, and  $a\in \WW Eb$ implies $\WW Ea=\WW Eb$ and $w(a,E)=w(b,E)$ for each $a,b\in V$
and $E\subs E^*$,
\item if $w(a,E)<w(a,E')$ for some $E\subs E'\subs E^*$ and $a\in V$, then 
$(E'\setm E)\cap [\WW Ea,V]\ne \empt$.
\end{enumerate}
\end{enumerate}
      \end{definition}

      \begin{theorem}\label{tm:braided-w}
            If $H=\<V,E^*\>$ is an infinite braided graph,
            $n\in {\omega}$, and   $w:V\times \mc P(E^*)\to n+1$ is a 
            non-trivial,  degree restricted, finitely determined and bounded function, then 
            the monotone graph property $R_n$
            \begin{displaymath}
            \text{``\em $w(a,E)=n$ for each $a\in V$''}
            \end{displaymath} 
            is not elusive in $H$.
            \end{theorem}

            \begin{proof}[Proof of Theorem \ref{tm:gen} from Theorem \ref{tm:braided-w}]
            If $P=D_m$ for some $m\ge 1$, define the function 
            $w:V\times \mc P(E^*)\to m+1$ as follows: 
            \begin{displaymath}
                  w(a,E)=\min(\ddeg_E(a),m). 
            \end{displaymath}
            Taking $K=m$, $M=1$ and  $\WW Ea =\{a\}$
            we obtain that $w$ is non-trivial, degree-restricted, finitely determined, and bounded. 
            
            By Proposition \ref{pr:braided} the graph $H$ is braided. 
            Thus, property (P) 
            \begin{displaymath}
                  \text{``\em $w(a,E)=m$ for each $a\in V$''}
                  \end{displaymath} 
                  is not elusive in $H$ by Theorem \ref{tm:braided-w}.
            But $P$ holds iff $D_m$ holds. Thus, we proved the theorem.

            If $P=C_n$ for some $n\ge 2$, 
            define the function 
            $w:V\times \mc P(E^*)\to n$ as follows: 
            \begin{displaymath}
                  w(a,E)=\min(\cc{a}{E}-1,n-1). 
            \end{displaymath}
            Taking $K=M=n-1$    and 
            \begin{displaymath}
            {\WW Ea }=\left\{\begin{array}{ll}
            {\CC{a}{E}}&\text{if $\cc{a}{E}< n$,}\\\\
            {\{a\}}&\text{if $\cc{a}{E}\ge  n$},
            \end{array}\right.
            \end{displaymath}
            we obtain that $w$ is non-trivial, degree-restricted, finitely determined, and bounded. 
            
            By Proposition \ref{pr:braided}, the graph $H$ is braided. 
            Thus, the property (P)
            \begin{displaymath}
                  \text{``\em $w(a,E)=n-1$ for each $a\in V$''}
                  \end{displaymath} 
                  is not elusive in $H$ by Theorem \ref{tm:braided-w}.
                  But $P$ holds iff $C_n$ holds. So we proved Theorem \ref{tm:gen}.   
            \end{proof}

Before proving Theorem \ref{tm:braided-w} we need some preparation. 

\begin{proposition}\label{pr:w-prop}
      Assume that $V$ is an infinite set, $H=\<V,E^*\>$ is a graph,  $n\in {\omega}$, and   
      $w:V\times \mc P(E^*)\to n+1$ is non-trivial, finitely determined and bounded function.     
 Then $w$ is
 \begin{enumerate}[(w1)]\addtocounter{enumi}{4}
      \item \label{w:monotone}  {\em monotone}, i.e.  
      $E_0\subs E_1$ implies $w(a,E_0)\le w(a,E_1)$ for $E_0\subs E_1\subs E^*$ and $a\in V$;
   \item \label{conti} {\em continuous},  i.e. if 
    $\<E_{\alpha}:{\alpha}<{\mu}\>\subs \mc P(E^*)$ is a $\subs$-increasing sequence, then 
   for each $a\in V$, we have 
   \begin{displaymath}
   w(a,\bigcup_{{\alpha}<{\mu}}E_{\alpha})=\max_{{\alpha}<{\mu}}w(a,E_{\alpha}).
   \end{displaymath} 
\item \label{stable} {\em stable} i.e. if $E\subs E^*$ and  $w(a,E)=w(b,E)=n$ for some 
$e=\{a,b\}\in E^*$, then 
$w(c,E\cup\{e\})=w(c,E)$ for each $c\in V$.  
 \end{enumerate}
\end{proposition}
 
\begin{proof}
If $E_0\subs E_1$, then for each $a\in V$,
\begin{multline*}
w(a,E_0)= \max \{w(a,E):E\in {[E_0]}^{<{\omega}}\}\le\\
\max \{w(a,E):E\in {[E_1]}^{<{\omega}}\} =w(a,E_1)
\end{multline*}
by (w\ref{w:finitely determined}). So (w\ref{w:monotone}) holds.

\medskip
If $\<E_{\alpha}:{\alpha}<{\mu}\>\subs \mc P(E^*)$ is a $\subs$-increasing sequence, then 
for each $a\in V$,
\begin{multline*}
      w(a,\bigcup_{{\alpha}<{\mu}}E_{\alpha})=\max\{w(a,E'):
      E\in {[\bigcup_{{\alpha}<{\mu}}E_{\alpha}]}^{<{\omega}}\}=\\
      \max\{\max\{w(a,E'):E'\in {[E_{\alpha}]}^{<{\omega}}\}:{\alpha}<{\mu}\}=
      \max_{{\alpha}<{\mu}}w(a,E_{\alpha}),
\end{multline*}
where the first and the third  equality hold by  
 (w\ref{w:finitely determined}), and the second holds because the sequence
 $\<E_{\alpha}:{\alpha}<{\mu}\>$ is increasing.
Thus, (w\ref{conti}) holds.       

\medskip
To prove (w\ref{stable})
assume  on the contrary that  $E\subs E^*$, $w(a,E)=w(b,E)=n$ for some $e=\{a,b\}\in E^*$, 
and  $w(c,E\cup\{e\})>w(c,E)$.
By (w\ref{w:bounded})(ii), $a\in \WW Ec $ or $b\in \WW Ec $. 
Thus, $w(c,E)=w(a,E)$ or $w(c,E)=w(b,E)$ by (w\ref{w:bounded})(i),
i.e. $w(c,E)=n$. But $w(c,E\cup\{e\})\le n$.  So $w(c,E\cup\{e\})>w(c,E)$ is not possible, 
(w\ref{stable}) holds.
\end{proof}

\begin{proof}[Proof of Theorem \ref{tm:braided-w}]

      Write ${\kappa}=|V|$.
Fix $M,K\in {\omega}$ and $W:V\times \mc P(E^*)\to {[V]}^{\le M}$       
as in Definition \ref{df:w}.

We will give a  winning strategy  for Alice in the game 
      $\mbb E_{H,R_n}$.
We divide the game into stages.

First Alice picks a finite covering set $L\subs V%
$.

      \noindent {\bf       Stage 1. }
      This stage is divided into $n$ substages. 

      Before the $i^{th}$ substage, 
      the players played $ m_i$ turns and determined a  pregraph $G_{m_i}$,
      such that  $m_i<{\kappa}$ and $w({\ell},{E_{m_i}})\ge i$ for each $\ell\in L$.
      
      Observe that the choice $m_0=0$ works because  
      $w ({\ell},\empt)=0\ge 0$  holds for each $\ell\in L$.

      \noindent {\bf Substage $i$. }

      {\em  
      Alice enumerates the edges of $H$ which contains at least one endpoint 
      from the finite set $\bigcup\{\WW{E_{m_i}}\ell:\ell\in L\}$
      as $\{e'_{\alpha}:{\alpha}<{\kappa}\}$ and in 
      the $m_i\dotplus{\alpha}$th turn Alice plans to play the undetermined pair $e'_{\alpha}$.
      }

If for each ${\alpha}<{\kappa}$ we have   
$w ({\ell},{G_{m_i+{\alpha}}})\le i$ for some $\ell\in L$, then 
$w( {\ell},{E_{m_i+{\kappa}}})\le i$ for some $\ell\in L$ because $w$ is monotone and continuous.

Since $[\WW{E_{m_i+{\kappa}}}{\ell},V]\subs P_{{m_i+{\kappa}}}$, we have 
$[\WW{E_{m_i+{\kappa}}}{\ell},V]\cap U_{{m_i+{\kappa}}}=\empt$,
and so $$\ww {\ell}{E_{m_i+{\kappa}}\cup U_{m_i+{\kappa}}}=\ww {\ell}{E_{m_i+{\kappa}}}\le i<n$$
by (w\ref{w:bounded})(ii).
Since $\gmax {G_{m_i+{\kappa}}}=\<V,E_{m_i+{\kappa}}\cup U_{m_i+{\kappa}}\>$,
it follows that $\gmax {G_{m_i+{\kappa}}}$ 
does not have 
property $R_n$.
Thus, the game terminates in at most $m_i+{\kappa}$ steps, and so 
Alice wins because $G_{m_i+{\kappa}}$ has undetermined edges. 

So we can assume that for some ${\alpha}<{\kappa}$ 
we have  $w( {\ell},{E_{m_i+{\alpha}}})\ge i+1$ for each $\ell\in L$, 
and after that step we declare that  the $i^{\rm th}$ substage is terminated, and 
$m_{i+1}=m_i\dotplus {\alpha}$.

      \bigskip

At the end of Stage 1, when we are after  $n$ many substages,
write ${\sigma}=m_{n}$ and observe that  
we have a pregraph $G_{{\sigma}}$  with $|P_{G_{{\sigma}}}|<{\kappa}$ such that 
$\ww {\ell}{E_{{\sigma}}}\ge n$ for each $\ell \in L$.  

For ${\alpha}<{\sigma}$ let  $e_{\alpha}$  denote the pair  Alice selected in the ${\alpha}$th turn.

\noindent {\bf       Stage 2.}
Let $N=M\cdot(K+|L|)+1$.
Alice picks  
$N+1$ pairwise disjoint  finite covering sets $\mc L=\{L_0,\dots, L_{N}\}$ from ${\kappa}\setm  
\bigcup\{e_{\alpha}:{\alpha}<{\sigma}\}$.
Let 
$$X=\bigcup \mc L
\text{ and } 
A={\kappa}\setm X.$$
Next
{\em 
Alice asks  all the undetermined pairs from   
$E^*\cap [A]^2$ in ${\kappa}$ turns. }

Since ${\sigma}\dotplus {\kappa}={\kappa}$ and $\www$ is monotone, so 
\begin{displaymath}
\text{$P_{{\kappa}}=E^*\cap [A]^2$ and $\ww v{E_{{\kappa}}}= n$ for each $v\in L\subs A$.
}
\end{displaymath}

Let 
\begin{displaymath}
B=\{{\beta}\in A: \ww {\beta}{E_{\kappa}}= n\} \text{ and } 
C=\{{\gamma}\in A: \ww {\gamma}{E_{\kappa}}<n\}.
\end{displaymath}
We have  $L\subs B$, and so   $B\ne \empt$.

\noindent {\bf       Stage 3. }

Alice  should distinguish two cases.

\noindent {\bf 
Case 3.1.}  $|B|={\kappa}$. %$|B|={\kappa}$ is infinite.  

Since $L_N$ is a covering set and $|B|={\kappa}$, Alice can pick $x\in L_N$ such that $|[\{x\},B]\cap E^*|={\kappa}$.
Let $Y=X\setm \{x\}$. 

\begin{figure}[ht]
\begin{tikzpicture}
 \path[draw,rounded corners] (0,0.5)  rectangle (3,1.5)  (3.2,1.5)  rectangle (6,0.5);  
 \path[draw,rounded corners,dotted] (-0.1,0.4)  rectangle  (6.1,1.6);  
\draw (3,0.4)  node[below right] {$A$} ;
 \path[draw]   (0,1.5) node[above left]  {$B$};
\path[draw]   (6,1.5) node[above right]  {$C$};
\path[draw] (0.5,0.8)  rectangle (1.5,1.2)    ;
\path[draw]  (1,1)  node {$L$};  
   \path[draw]  (1.5,1.5)  node[above] {$\ww v{E_{\kappa}}= n$};
 \path[draw]  (4.5,1.5)  node[above] {$\ww v{E_{\kappa}}< n$};
 \path[draw] (7,0.8)  rectangle (8,1.2)    ;
 \path[draw] (8.5,0.8)  rectangle (9.5,1.2)    ;
 \path[draw] (10,0.8)  rectangle (11,1.2)    ;
 \path[draw]  (7.5,1)  node {$L_0$};  
 \path[draw]  (10.5,1)  node {$L_{N-1}$};  
 \path[draw,fill]  (11.5,1) circle (2pt) node[above] {$x$};
 \path[draw,dotted] (6.8,0.6)  rectangle (11.2,1.4)  node[above]  {$Y$};  
\end{tikzpicture}
\caption{}\label{fig:c2x}
\end{figure}

\noindent {\bf       Substage 3.1.1 }

{\em First Alice asks all the undetermined pairs 
from $[Y]^2$ in finitely many steps,  then the undetermined pairs
$[C,Y]$ in at most ${\kappa}$ turns. } ($C$ can be the empty set).

So for some ${\rho}\le {\kappa}$
\begin{displaymath}
P_{{{\kappa}\dotplus {\rho}}}=
E^*\cap ({[A]}^{2}\cup {[C\cup Y]}^{2})= E^*\cap ({[{\kappa}\setm \{x\}]}^{2}\setm [B,Y]).
\end{displaymath}

\begin{claim}
If $\ww y{E_{{\kappa}\dotplus {\rho}}}= n$
for some $y\in Y$, then Alice can win.
\end{claim}

\begin{proof}
Since $L$ is a covering set, Alice can pick  $\ell \in L$ with $\{y,\ell\}\in E^*$.

{\em Alice enumerates the undetermined pairs of $G_{{\kappa}\dotplus {\rho}}$  with 
the exception of 
$\{y,\ell\}$ 
as 
$\{e'_{\gamma}:{\gamma}<{\kappa}\}$,
and 
she asks the pair $e'_{\gamma}$ in the turn ${\kappa}\dotplus {\rho}\dotplus {\gamma}$.}

Writing  ${\delta}={\kappa}\dotplus {\rho}\dotplus {\kappa}$ 
we have 
\begin{displaymath}
      \text{      
      $U_{{{\delta}}}=\{\{\ell,y\}\}$ and 
      $\ww y{E_{{\delta}}}= 
      \ww \ell{E_{{\delta}}}= n$. 
      } 
      \end{displaymath}
      So $\gmin G_{{\delta}}=\<V, E_{{{\delta}}}\>$ has property $R_n$ iff 
      $\gmax {G_{\delta}}=\<V, E_{{{\delta}}}\cup\{\{\ell,y\}\}\>$ has property $R_n$
      by (w\ref{stable}). Thus, the game finishes after at most  ${\delta}$ turns and $U_{{\delta}}\ne \empt$.
        So  Alice wins. 
\end{proof}

So we can assume that 
\begin{displaymath}
\text{$\ww y{E_{{\kappa}\dotplus {\rho}}}<n$
for each $y\in Y$. 
}
\end{displaymath}

\noindent {\bf       Substage 3.1.2 }

%Since $C$ is finite, $|B|={\kappa}$.
{\em Alice enumerates the pairs  $[B,Y]$  in type ${\kappa}$ as $\{e'_{\zeta}:{\zeta}<{\kappa}\}$, 
and in the turn ${\kappa}\dotplus {\rho}\dotplus  {\zeta}$ she plans to play $e'_{\zeta}$. }

For each $i=0,\dots, N-1$ pick $\ell_i\in L_i$ such that 
$|[\{\ell_i\},B]\cap E^*|={\kappa}$.

\begin{claim}
If for some ${\kappa}\dotplus {\rho}\le {\eta}<{\kappa}\dotplus {\rho}\dotplus {\kappa}$  
we have that $\ww {\ell_i}{E_{\eta}}= n$
for some $0\le i< N$, then Alice can win.
\end{claim}

\begin{proof}
After realizing that $\ww {\ell_i}{E_{\eta}}= n$ after the turn ${\eta}$, Alice changes her plan 
how to play after that turn. 
      Since $|[\{\ell_i\},B]\cap E^*|={\kappa}$, we
can pick a pair $e=\{\ell_i,b\}\in [Y,B]$ such that 
$e$ is undetermined in $G_{\eta}$. 
{\em In the next ${\kappa}$ turns  Alice asks all the undetermined pairs of $G_{\eta}$ with the 
exception of $e$}.
Then 
\begin{displaymath}
\text{      
$U_{{\eta}\dotplus{\kappa}}=\{e\}$ and 
$\ww {\ell_i}{G_{{\eta}\dotplus{\kappa}}}= n$ and $\ww b{G_{{\eta}\dotplus{\kappa}}}= n$. 
} 
\end{displaymath}
So $\gmin G_{{\eta}\dotplus {\kappa}}=\<V, E_{{{\eta}\dotplus{\kappa}}}\>$ has property $R_n$ iff 
$\gmax {G_{{\eta}\dotplus {\kappa}}}=\<V, E_{{{\eta}\dotplus{\kappa}}}\cup\{e_{\gamma}\}\>$ has 
property $R_n$  by  (w\ref{stable}). 
Thus, the game finishes after at most  ${\eta}+{\kappa}$ turns and $U_{{\eta}+{\kappa}}\ne \empt$.
Thus,  Alice wins. 
\end{proof}

So we can assume that Bob answered in such a way that after 
${\eta}={\kappa}\dotplus {\rho}\dotplus{\kappa}$ turns 
\begin{displaymath}
\text{$\ww {\ell_i}{E_\eta}<n$ for each  $i< N$, }
\end{displaymath}
and 
\begin{displaymath}%\tag{$*$}
P_{G_{\eta}}=E^*\cap {[A]}^{2}\cup {[Y]}^{2}\cup [A,Y]=E^*\cap {[A\cup Y]}^{2}.
\end{displaymath} 
Let $$Z=\{\ell_i: i< N\}.$$

\bigskip

      \noindent {\bf 
      Case 3.2.}   $|B|<{\kappa}$. 
      
      Then   $|C|={\kappa}$.  
Pick  an arbitrary $x\in L_N$. 
      Let $Y=X\setm \{x\}$.

\noindent {\bf Substage 3.2.1}

{\em First Alice enumerates  all the pairs from 
 $[Y]^2$ and  the undetermined pairs in  $[C,Y]$ 
as $\{e'_{\gamma}:{\gamma}<{\kappa}\}$, and she asks the pair
$e'_{\gamma}$ in the turn ${\kappa}\dotplus {\gamma}$.} 
Consider the set 
 $$C'=\{y\in C: \ww y{E_{{\kappa}\dotplus {\kappa}}}<n \}.$$

 \begin{figure}[ht]
      \begin{tikzpicture}
       \path[draw,rounded corners] (0,0.5)  rectangle (3,1.5)  (3.2,1.5)  rectangle (8,0.5);
       \path[draw,rounded corners]  (3.3,1.3) rectangle (6,0.7);  
       \path[draw,rounded corners,dotted] (-0.1,0.4)  rectangle  (8.1,1.6);  
      \draw (3,0.4)  node[below right] {$A$} ;
       \path[draw]   (0,1.5) node[above left]  {$B$};
      \path[draw]   (8,1.5) node[above right]  {$C$};
      \path[draw]   (6,1) node[below right]  {$C'$};
      \path[draw] (0.5,0.8)  rectangle (1.5,1.2)    ;
      \path[draw]  (1.5,1)  node[above right] {$L$};  
      \path[draw,fill]  (1,1) circle (2pt) node[above] {$\ell$};
      \path[draw]  (1.5,1.5)  node[above] {$\ww v{E_{\kappa}}= n$};
         \path[draw]  (4.5,0.7)  node[above] {$\ww v{E_{{\kappa}+{\kappa}}}< n$};
       \path[draw]  (4.5,1.5)  node[above] {$\ww v{E_{\kappa}}< n$};
       \path[draw,fill]  (11.5,1) circle (2pt) node[above] {$x$};
       \path[draw,fill]  (9.5,1) circle (2pt) node[above] {$y$};
       \path[draw,dotted] (8.8,0.6)  rectangle (11.2,1.4)  node[above]  {$Y$};  
      \end{tikzpicture}
      \caption{}\label{fig:c2xx}
      \end{figure}

\begin{claim}
If $|C'|<{\kappa}$, then Alice can win.
\end{claim}

\begin{proof}
 If $|C'|<{\kappa}$, then 
 for each $c\in C\setm C'$ we have 
 $\ww{c}{E_{{\kappa}+{\kappa}}}>\ww{c}{E_{\kappa}} $, so 
 by (w\ref{w:bounded})(ii) we can find
 an edge $\{y_c,a_c\}\in E_{{\kappa}+{\kappa}}\setm E_{\kappa}$ such that
 $a_c\in \WW{E_{\kappa}}c$.
 
Since $|C\setm C'|={\kappa}$, by (w\ref{w:bounded})(i) we can find 
$C''\in {[C\setm C']}^{{\kappa}}$ such that $\{\WW{E_{\kappa}}c:c\in C''\}$
are pairwise disjoint  and $Y\cap \bigcup\{\WW{E_{\kappa}}c:c\in C''\}=\empt$.

 Then $y_c\in Y$ for $c\in C''$. 
Since $Y$ is finite, there is $y\in Y$ such that 
$|\{c\in C'':y_c=y\}|={\kappa}$ and so 
$\ddeg_{E_{{\kappa}\dotplus {\kappa}}}(y)={\kappa}$. 
Thus,  $\ww y{E_{{\kappa}\dotplus {\kappa}}}= n$ by (w\ref{w:deg}). 

Pick $\ell\in L$ such that $\{\ell,y\}\in E^*$.
 
{\em Next
Alice enumerates  all undetermined pairs  of $G_{{\kappa}\dotplus {\kappa}}$ with the exception of  
$\{\ell,y\}$ as  $\{e_{\gamma}:{\gamma}<{\kappa}\}$,
and she asks $e_{\gamma}$ in the turn ${\kappa}\dotplus {\kappa}\dotplus {\gamma}$.}

Let ${\eta}={\kappa}\dotplus {\kappa}\dotplus {\kappa}$.

Then 
\begin{displaymath}
      \text{      
      $U_{{\eta}}=\{\{\ell,y\}\}$ and 
      $\ww y{E_{{\eta}}}= n$ and 
      $\ww \ell{E_{{\eta}}}= n$. 
      } 
      \end{displaymath}
      So $\gmin G_{\eta}= \<V, E_{{{\eta}}}\>$ has property $C_n$ iff 
      $\gmax {G_{\eta}}=\<V, E_{{{\eta}}}\cup\{e_{\gamma}\}\>$ has property $R_n$
      by  (w\ref{stable}).
      Thus, the game terminates after at most  ${\eta}$ turns and $U_{{\eta}}\ne \empt$.
      Thus,  Alice wins. 
\end{proof}

So we can assume that $|C'|={\kappa}$. 

\noindent {\bf Substage 3.2.2}

{\em Alice asks all the undetermined pairs from $[B,Y]$ in the next ${\sigma}\le {\kappa}$
turns.} 

If 
$\ww c{E_{{\kappa}+{\kappa}+{\sigma}}}>\ww c{E_{{\kappa}+{\kappa}}}$ for some $c\in C'$
then $\WW {E_{{\kappa}+{\kappa}}}c\cap B= \empt$ by (w\ref{w:bounded})(i) and so 
$\WW {E_{{\kappa}+{\kappa}}}c\cap Y\ne \empt$ by (w\ref{w:bounded})(ii) 
and so $c\in \WW{E_{{\kappa}+{\kappa}}}{y}$ for some $y\in Y$ by (w\ref{w:bounded})(i).
 So 
\begin{displaymath}
      |\{c\in C':\ww c{E_{{\kappa}+{\kappa}+{\sigma}}}>
      \ww c{E_{{\kappa}+{\kappa}}}\}|\le |Y|\cdot M <{\omega}.
\end{displaymath}
 i.e. writing ${\eta}={\kappa}\dotplus{\kappa}\dotplus{\sigma}$
we have that there is an (infinite) $C''\subs C'$ with $|C'\setm C''|\le |Y|\cdot M)$ 
such that   
\begin{displaymath}
\forall v\in C''\quad \ww v{E_{\eta}}=\ww v{E_{{\kappa}+{\kappa}}}<n.
\end{displaymath} 
Fix $Z\in {[C'']}^{N}$.
   
\bigskip
 \noindent {\bf Stage 4.}

 Both in Case 3.1 and in Case 3.2 
 after ${\eta}$ turns we have 
 \begin{displaymath}
 P_{{\eta}}=E^*\cap  {[{\kappa}\setm \{x\}]}^{2}
 \end{displaymath}
 and 
 there is a set $Z\in {[{\kappa}\setm \{x\}]}^{N}$
 such that $\ww z{E_{\eta}}<n$ for each $z\in Z$.

 Since $|Z|>M(K+|L|)$,
 applying  the fact that  for each $\{z,z'\}\in {[V]}^{2}$ either $\WW{E_{\eta}}{z}=\WW{E_{\eta}}{z'}$
or  $\WW{E_{\eta}}{z}\cap\WW{E_{\eta}}{z'}=\empt$  by (w\ref{w:bounded})(i),
we can pick elements %$z_0,\dots z_{K+|L|-1},z_{K+|L|}$  
$\{z_i:i\le K+|L|\}$  
from  $Z$ such that 
$\{\WW{E_{\eta}}{z_i}: i\le K+|L|\}$ are pairwise disjoint. 
Write $Z_i=W_{E_{\eta}}(z_i)$.
We can assume that 
\begin{displaymath}
      \tag{x}(L\cup \{x\})\cap  \bigcup_{i<K}Z_i=\empt.
\end{displaymath}

Pick $\ell \in L$ such that $\{x,\ell\}\in E^*$.

 \begin{figure}[ht]
\begin{tikzpicture}
\path[draw] (0,0) rectangle (10,3) ;      
\path[draw,fill]  (11,1.5)  circle  (1.5pt)  node[above] {$x$} ;
\path[draw,fill]  (1,1.5)  circle  (1.5pt)  node[above] {$\ell$} node[below right] 
{$\ww \ell{E_{\eta}}\ge n$} ;
\path[draw] (4,1) rectangle (9,2)  node[above right] {$Z$};      
\path[draw] (4.5,1.0) circle  (11.5pt)  node[above] {$Z_0$} ;
\path[draw] (5.5,1.0) circle  (11.5pt)  node[above] {$Z_2$} ;
\path[draw] (6.5,1.0)   node[above] {$\dots$} ;
\path[draw] (7.5,1.0) circle  (11.5pt)  node[above] {$Z_{K-1}$} ;
\path[draw]  (6,2.5) node[below]  {$\ww v{E_{\eta}}<n$}; 
\end{tikzpicture}
\caption{}\label{fig:c3}
\end{figure}

Alice enumerates the undetermined pairs 
of $G_{\eta}$ with the exception of $\{x,\ell\}$
as  $\{e_{\gamma}:{\gamma}<{\delta}\}$,
and she asks $e_{\gamma}$ in the turn ${\eta}+{\gamma}$. 

If $[Z_i,\{x\}]\cap E_{\eta+{\delta}}=\empt$,
then  $[Z_i,\{x\}]\cap (E_{{{\eta}+{\delta}}}\cup\{\{x,\ell\}\})=\empt$ as well by (x), and so 
$$\ww z{{E _{{\eta}+{\delta}}}\cup\{\{x,\ell\}\}}=\ww z{{E _{{\eta}+{\delta}}}}=\ww z{E_{\eta}}<n$$
for each $z\in Z_i$
by (w\ref{w:bounded})(ii).
Since 
\begin{displaymath}
\text{$U_{{\eta}+{\delta}}=\{\{x,l\}\}$ and $\gmax {\<V,E_{{\eta}+{\delta}}\>}=
\<V,{{E _{{\eta}+{\delta}}}\cup\{x,\ell\}} \>$ does not have property $R_n$},
\end{displaymath}
so the game terminates at least ${\eta}+{\delta}$ steps and $U_{{\eta}+{\delta}}\ne \empt$,
so Alice wins.

So we can assume that $[Z_i,\{x\}]\cap E_{{{\eta}+{\delta}}}\ne \empt$ for each $i<K$.
So  
$\ddeg_{G_{{\eta}+{\delta}}}(x)\ge K$, and so $\ww{x}{E_{{\eta}+{\delta}}}= n$ 
by (w\ref{w:deg}).
Since $U_{{{\eta}+{\delta}}}=\{\{\ell,x\}\}$ and 
$\ww{\ell}{E_{{\eta}+{\delta}}}  = n$, 
\begin{multline*}
\text{       
$\gmin G_{\eta+{\delta}}= \<{\kappa},E_{{\eta}+{\delta}}\>$ has property $R_n$ iff}\\
\text{  
$\gmax {G_{{\eta+{\delta}}}}=\<{\kappa},E_{{{\eta}+{\delta}}}\cup \{\{\ell,x\}\}\>$ 
has property $R_n$}
\end{multline*}
by  (w\ref{stable}). Thus, the game terminates after at most  ${\eta}+{\delta}$ turns and 
$U_{{\eta}+{\delta}}\ne \empt$.
Thus, Alice wins. This completes the proof of Theorem \ref{tm:braided-w}.
\end{proof}

\noindent {\bf Problems:}
 (1)
Are the following properties (strongly) elusive for infinite vertex sets:
\begin{enumerate}[(i)]
\item ``$G$  contains $P_3$'',
\item ``$G$  contains $K_n$'' for some $3\le n< {\omega}$,
\item ``$G$ is  not bipartite'',
\item ``$G$ is connected.''
\end{enumerate}
(2) Is it true that a property is elusive iff it is strongly elusive? 
Is it true under the Axiom of Determinacy?

\noindent (3) Find a reasonable conjecture to characterize (large classes of) monotone graph properties which 
are elusive for infinite vertex sets.

\section{On the infinite version of the Aanderaa-Rosenberg conjecture}
\label{sc:weak-elusive}

We have seen in the previous section that Alice does not have to ask all
the pairs in an infinite vertex set to decide if the hidden graph $G$ has property 
$D_n$ or $C_m$.  

In contrast to these results,
in this section we show that Alice, the seeker should check "lots of edges" 
before she can decide whether a graph has certain monotone graph properties.

\begin{definition}\label{df:w-elusive1-tartalek} 

      Let $R$ be a graph property,  $H=\<V,E^*\>$ be a graph, and 
      $\mc F\subs \mc P(E^*)$ be upward closed. We will say that 
      $\mc F$ is the {\em family of large edge sets}. 
%      The players will construct a subgraph of $H=\<V,E^*\>$.

      Define the game $\mbb E_{H,R,\mc F}$ between two players, Alice and Bob, as follows: 

      The gameplay is the same as in the game $\mbb E_{H,R}$, but we modify the rule which determines 
      when Bob wins: so we keep
      rules (1)-(4) from Definition \ref{df:w-elusive1}, but we replace (5) with 
      \begin{enumerate}[(1')]\addtocounter{enumi}{4}
            \item 
Bob wins iff $P_{\beta}\in \mc F$.      
(Informally, Bob wins if he can force Alice to ask a large set of edges.) 
\end{enumerate}

If we do not have any restriction on edges, i.e. if
$H=K_V=\bigl\langle V,{[V]}^{2}\bigr\rangle$, 
 then we write $\mbb E_{V,R,\mc F}$ for $\mbb E_{K_V,R,\mc F}$. 
     
Observe that $\mbb E_{V,R}=\mbb E_{V,R,\{{[V]}^{2}\}}$.

      \end{definition}

\begin{definition}
      Given a   vertex set $V$ and  a family $\mc F\subs \mc P({[V]}^{2})$,  
      a graph property $R$ is said to be {\em $\mc F$-hard on $V$}  iff Alice does not 
      have a winning strategy in the game $\mbb E_{V,R,\mc F}$.
      \end{definition}

As we recalled earlier, 
Rivest and Vuillemin \cite{RiVu}   proved that if $R$ is a non-trivial monotone 
      graph property, then  for each finite vertex set $V$ the seeker should ask $\Omega(|V|^2)$ 
      many edges in the worst case. 
We think that this result yields naturally the following problem for infinite sets:

\begin{definition} 
      Let 
      \begin{displaymath}
      \isf=\{F\subs {[{\omega}]}^{2}: \exists W\in {[{\omega}]}^{{\omega}}\ {[W]}^{2}\subs F\}.
      \end{displaymath}
      \end{definition}

      \begin{problem}
            Assume that $P$ is a non-trivial, monotone graph property.
            Is it true that $P$   is $\isf$-hard on ${\omega}$?
            \end{problem}

       \begin{definition}\label{df:scorpion}
       A {\em scorpion graph on a vertex set $V$} is a graph $G=\<V,E\>$
       such that it contains three special vertices, 
       the {\em sting} denoted by $s_G$,
       the {\em tail} denoted by $t_G$, and the {\em body} 
       denoted by $b_G$ such that 
       \begin{enumerate}[(a)]
       \item 
       the sting is connected only to the tail, 
       \item the tail is connected only
       to the sting and to the body, and 
       \item
       the body is connected to all vertices except the sting. 
      \end{enumerate}
%       \begin{figure}[ht]
%             \begin{tikzpicture}
% \path[draw,fill] (1,1) circle (2pt) node[above] {sting};
% \path[draw,fill] (3,1) circle (2pt) node[above] {tail};
% \path[draw,fill] (5,1) circle (2pt) node[above] {body};
% \path[draw]  (1,1) -- (3,1) -- (5,1);
% \path[draw,rounded corners,dotted]  (6,-1) rectangle (8,3) ;
% \foreach \i in {-0.5,0.5,...,2.5}
% {  \path[draw,fill] (5,1) -- (7,\i) circle (2pt) ;  
% };
% \path[draw] (7,-0.5) .. controls (7.5,0.5) .. (7,1.5);     
% \path[draw] (7,1.5) .. controls (7.5,2) .. (7,2.5);     
% \end{tikzpicture}
%             \caption{}
%        \end{figure}     
       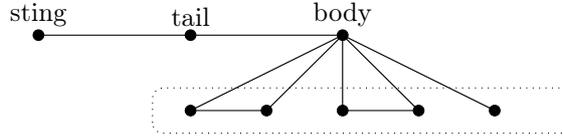
\begin{figure}[ht]
            \begin{tikzpicture}
\path[draw,fill] (1,1) circle (2pt) node[above] {sting};
\path[draw,fill] (3,1) circle (2pt) node[above] {tail};
\path[draw,fill] (5,1) circle (2pt) node[above] {body};
\path[draw]  (1,1) -- (3,1) -- (5,1);
\path[draw,rounded corners,dotted]  (2.5,-0.3) rectangle (8,0.3) ;
\foreach \i in {3,4,5,6,7}
{  \path[draw,fill] (5,1) -- (\i,0) circle (2pt) ;  
};
\path[draw] (3,-0) -- (4,0)  (5,-0) -- (6,0) ; 
%\path[draw] (4,0) .. controls (5.5,-0.5) .. (7,0);     
\end{tikzpicture}
            \caption{Scorpion graph}
       \end{figure}     
%The sting  of a scorpion graph $G$ is denote by $s_G$, the tail  by 
%$t_G$, and the body  by $b_G$.
      \end{definition}     

In \cite{LBB} it was proved that  that 
if $V$ is a finite set, then 
the 
"{\em being a scorpion graph on $V$}"
property can be decided by the seeker in $O(|V|)$ steps. 

\newcommand{\sid}[1]{\mathcal J_{#1}}

\newcommand{\infi}[1]{\Omega(#1)}

\begin{definition}
If $E\subs {[{\omega}]}^{2},$ let 
\begin{displaymath}
\infi{E}=\{v\in {\omega}:\deg_E(v)={\omega}\}, 
\end{displaymath}      
and for $n\in {\omega}$ let 
\begin{displaymath}
\sid{n}=\{E\subs {[{\omega}]}^{2}: |\infi{E}|\ge n\}.
\end{displaymath}
\end{definition}

\begin{theorem}\label{tm:scorpion}
The property $(S)$
\begin{displaymath}
      \text{"{\em being a scorpion graph"}}
\end{displaymath}
is not $\sid 5$-hard on the vertex set ${\omega}$.
\end{theorem}

Before proving the theorem we need some preparation. 

\newcommand{\BBB}[1]{B_{#1}}
\newcommand{\SSS}[1]{S_{#1}}
\newcommand{\TTT}[1]{T_{#1}}

\begin{definition}
If $G=\<{\omega},E,N\>$ is a pregraph, let 
\begin{align*}
\BBB{G}=&\{v\in {\omega}: \deg_N(v)\le 1\},\\
\SSS{G}=&\{v\in {\omega}: \deg_E(v)\le 1\},\\
\TTT{G}=&\{v\in {\omega}: \deg_E(v)\le 2\}.\\
\end{align*}
\end{definition}

The following two lemmas are straightforward from the definition.
\begin{lemma}\label{lm:extension}
If $G=\<{\omega},E,N\>$ is a pregraph, and $G\le H$ is a scorpion graph,
then $b_H\in \BBB{G}$, $s_H\in \SSS{G}$ and $t_H\in \TTT{G}$.
\end{lemma}

\begin{lemma}\label{lm:s-char}
If $G=\<{\omega},E,N\>$ is a pregraph, $s',t',b'\in {\omega}$, 
$[\{s',t',b'\},{\omega}]\subs E\cup N$, then the following statements are equivalent:
\begin{enumerate}[(1)]
\item $G'=\<{\omega},E\>$ is a scorpion graph
with $s_{G'}=s'$,
$t_{G'}=t'$ and $b_{G'}=b'$,
\item some extension $H$ of $G$ is a scorpion graph with $s_{H}=s'$,
$t_{H}=t'$ and $b_{H}=b'$,
\item every extension $H$ of $G$ is a scorpion graph with $s_{H}=s'$,
$t_{H}=t'$ and $b_{H}=b'$.
\end{enumerate}
\end{lemma}

% Actually, we will prove a stronger result. 

%       \begin{definition}
%             For $n\in {\omega}$ let 
%             \begin{displaymath}
%             \isk n=\{F\subs {[{\omega}]}^{2}: 
%             \exists A\in {[{\omega}]}^{n}\ ({[A]}^{2}\subs F)
%             \}.
%             \end{displaymath}
%             \end{definition}
           
% Clearly 
% \begin{displaymath}
% \isf \subs \bigcap_{n\in {\omega}}\isk n.
% \end{displaymath}

% The following theorem clearly implies Theorem \ref{tm:scorpion}. 

% \begin{theorem}\label{tm:scorpionpluss}
%       The property (S)
%       \begin{displaymath}
%             \text{"{\em being a scorpion graph"}}
%       \end{displaymath}
%       is not $\isk 6$-hard on the vertex set ${\omega}$.
%       \end{theorem}

\begin{proof}[Proof of theorem \ref{tm:scorpion}]
      
      We will use the following notation. If the players construct 
      the sequence $\<G_{\alpha}:{\alpha}\le {\beta}\>$ of pregraphs, 
      then we write $G_{\alpha}=\<{\omega},E_{\alpha},N_{\alpha}\>$,
      $P_{\alpha}=P_{G_{\alpha}}=E_{\alpha}\cup N_{\alpha}$, 
      $U_{\alpha}=U_{G_{\alpha}}={[{\omega}]}^{2}\setm P_{\alpha}$,  
      $\SSS{\alpha}=\SSS{G_{\alpha}}$, $\BBB{\alpha}=\BBB {G_{\alpha}}$
      and $\TTT{\alpha}=\TTT{G_{\alpha}}$.

      Let 
      \begin{displaymath}
      \mbb S=\{G\ : \ \text{$G$ is a scorpion graph on ${\omega}$}\}.
      \end{displaymath}
      
      Before  describing a winning strategy of Alice  in the game  
      $\mbb E_{S,{\omega},\sid 5}$, we need to prove the following lemma.

\begin{lemma}\label{lm:alg}
      Assume that Alice and Bob played so far ${\alpha}$ turn in the game $\mbb E_{S,{\omega},\sid 5}$,
      and  they constructed the pregraph $G_{\alpha}$.
      %obtained the sequence  $\<G_{\zeta}:{\zeta}\le {\alpha}\>$ 
      %of pregraphs.  
     
      If either 
 \begin{enumerate}[(1)]
 \item[(a)] there is $s'\in {\omega}$ such that  $s_G=s'$
 for each  $G\in \mbb S$ with  $G_{\alpha}\le G$,
 \end{enumerate}
 or    
 \begin{enumerate}[(1)]
      \item[(b)] there is $b'\in {\omega}$ such that $b_G=b'$
 for each  $G\in \mbb S$ with  $G_{\alpha}\le G$,
 \end{enumerate}
      then Alice can 
      play in such a way that the game terminates in ${\beta}\le {\alpha}+{\omega}+{\omega}+{\omega}$
      turns and $|\infi{P_{\beta}}\setm \infi{P_{\alpha}}|\le 3$.
\end{lemma}

\begin{proof}[Proof of Lemma \ref{lm:alg}]
Assume first that (a) holds: $s'\in {\omega}$ such that  $s_G=s'$
for each  $G\in \mbb S$ with  $G_{\alpha}\le G$. 

The strategy of Alice will be divided into three stages. 

\noindent{\bf Stage A.}

{\em Alice enumerates  the undetermined pairs 
$U_{\alpha}\cap [\{s'\},{\omega}]$
as $\{e'_\ell:\ell<M_s\}$ for some $M_s\le {\omega}$, and 
Alice plays $e'_\ell$ in step ${\alpha}+\ell$ for $\ell<M_s$.} 

\smallskip

Write ${\sigma}={\alpha}+M_s$.
Clearly $\infi{P_{\sigma}}\setm \infi{P_{\alpha}}\subs \{s'\}$.
    Let 
\begin{displaymath}
T=\{t\in {\omega}:\{t,s'\}\in E_{\sigma}\}.
\end{displaymath}

If $G\ge G_{\sigma}$ is a scorpion graph, then $s_G=s'$, and so 
$T=\{t_G\}$.

Thus, if $|T|\ne 1$, then there is no scorpion graph $G$ which extends $G_{\sigma}$.
So the game terminates after at most ${\sigma}$ turns and $\infi{P_{\sigma}}\setm \infi{P_{\alpha}}\subs \{s'\}$.

So we can assume that  $T=\{t'\}$, and  we know that 
\begin{displaymath}
\tag{B}\text{if $G \ge G_{{\sigma}}$ is a scorpion graph, then $s_G=s'$ and    $t_G=t'$.}
\end{displaymath}

\noindent{\bf Stage B.}

{\em Alice enumerates  the undetermined pairs 
$U_{{\sigma}}\cap [\{t'\},{\omega}]$
as $\{e''_\ell:\ell<M_t\}$ for some $M_t\le {\omega}$, and 
Alice plays $e''_\ell$ in step ${\sigma}+\ell$ for $\ell<M_t$.}

\medskip
Write ${\rho}={\sigma}+M_t$.
Clearly $\infi{P_{\rho}}\setm \infi{P_{\alpha}}\subs \{s',t'\}$.

Let
\begin{displaymath}
B=\{b\in {\omega}:\{b,t'\}\in E_{\rho}, b\ne s'\}.
\end{displaymath}

If $G\ge G_{\rho}$ is a scorpion graph, then  $s_G=s'$  and $t_G=t'$, and so 
$B=\{b_G\}$.

Thus, if $|B|\ne 1$, then there is no scorpion graph $G$ which extends $G_{\rho}$. 
So the game terminates after at most ${\rho}$ turns and 
$\infi{P_{\rho}}\setm \infi{P_{\alpha}}\subs \{s',t'\}$.

So we can assume that  $B=\{b'\}$, and  we know that 
\begin{displaymath}
\tag{C}\text{if $G \ge G_{\rho}$ is a scorpion graph, then $s_G=s'$, $t_G=t'$  and $b_G=b'$.}
\end{displaymath}

\noindent{\bf Stage C.}

{\em Alice enumerates  the undetermined pairs 
$U_{{\rho}}\cap [\{b'\},{\omega}]$
as $\{e^*_\ell:\ell<M_b\}$ for some $M_b\le {\omega}$, and 
Alice plays $e^*_\ell$ in step ${\rho}+\ell$ for $\ell<M_b$.} 

\medskip
Write ${\nu}={\rho}+M_b$.
Clearly $\infi{P_{\nu}}\setm \infi{P_{\alpha}}\subs  \{s',t',b'\}$.

Since (C) holds and  
$[\{s',t',b'\},{\omega}]\subs P_{\nu}$,
by Lemma \ref{lm:s-char} either 
every extension of $G_{\nu}$ is a scorpion graph, or  
no extension of $G_{\nu}$ is a scorpion graph.

So the game terminates after at most ${\nu}$ turns, 
${\nu}\le {\alpha}+{\omega}+{\omega}+{\omega}$
and $\infi{P_{\nu}}\setm \infi{P_{\alpha}}\subs  \{s',t',b'\}$.

\medskip

Assume now that (b) holds: $b'\in {\omega}$ such that $b_G=b'$
for each  $G\in \mbb S$ with  $G_{\alpha}\le G$. 

The strategy of Alice will be divided into stages.

\noindent {\bf Stage A.}
{\em Enumerate $U_{\alpha}\cap [\{b'\},{\omega}]$ as $\{e'_i:i<M_b\}$
for some $M_b\le {\omega}$ and in the next $M_b$ steps Alice asks 
$\{e'_i:i<M_b\}$. }

\medskip
Write ${\sigma}={\alpha}+M_b$. 
Clearly $\infi{P_{\sigma}}\setm \infi{P_{\alpha}}\subs \{b'\}$.
    Let 
\begin{displaymath}
S=\{s\in {\omega}:\{b',s'\}\in N_{\sigma}\}.
\end{displaymath}

If $G\ge G_{\sigma}$ is a scorpion graph, then $b_G=b'$, and so 
$S=\{s_G\}$.

Thus, if $|S|\ne 1$, then there is no scorpion graph $G$ which extends $G_{\sigma}$, and 
so the game terminates after at most ${\sigma}$ turns and 
$\infi{P_{\sigma}}\setm \infi{P_{\alpha}}\subs \{s'\}$.

So we can assume that  $S=\{s'\}$, and  we know that 
\begin{displaymath}
\tag{B''}\text{if $G \ge G_{\sigma}$ is a scorpion graph, then $b_G=b'$ and    $s_G=s'$.}
\end{displaymath}

\noindent{\bf Stage B.}

{\em Alice enumerates  the undetermined pairs 
$U_{{\sigma}}\cap [\{s'\},{\omega}]$
as $\{e''_\ell:\ell<M_s\}$ for some $M_s\le {\omega}$, and 
Alice plays $e''_\ell$ in step ${\sigma}+\ell$ for $\ell<M_s$.}

\medskip
Write ${\rho}={\sigma}+M_s$.
Clearly $\infi{P_{\rho}}\setm \infi{P_{\alpha}}\subs \{b',s'\}$.

Let
\begin{displaymath}
T=\{t\in {\omega}:\{t,s'\}\in E_{\rho}\}.
\end{displaymath}

If $G\ge G_{\rho}$ is a scorpion graph, then  $s_G=s'$   and so 
$T=\{t_G\}$.

Thus, if $|T|\ne 1$, then there is no scorpion graph $G$ which extends $G_{\rho}$, and  
so the game terminates after at most ${\rho}$ turns
and $\infi{P_{\rho}}\setm \infi{P_{\alpha}}\subs \{b',s'\}$.

So we can assume that  $T=\{t'\}$, and  we know that 
\begin{displaymath}
\tag{C''}\text{If $G \ge G_{\rho}$ is a scorpion graph, then $s_G=s'$, $t_G=t'$  and $b_G=b'$.}
\end{displaymath}
\noindent{\bf Stage C.}

{\em Alice enumerates  the undetermined pairs 
$U_{{\rho}}\cap [\{t'\},{\omega}]$
as $\{e^*_\ell:\ell<M_t\}$ for some $M_t\le {\omega}$, and 
Alice plays $e^*_\ell$ in step ${\rho}+\ell$ for $\ell<M_t$.}

\medskip

Write ${\nu}={\rho}+M_t$.
Clearly $\infi{P_{\nu}}\setm \infi{P_{\alpha}}\subs \{s',t',b'\}$.

Since (C'') holds and  
$[\{s',t',b'\},{\omega}]\subs P_{\nu}$,
by Lemma \ref{lm:s-char} either 
every extension of $G_{\nu}$ is a scorpion graph, or  
no extension of $G_{\nu}$ is a scorpion graph.

So the game terminates after at most ${\nu}$ turns,
${\nu}\le {\alpha}+{\omega}+{\omega}+{\omega}$ and 
$\infi{P_{\nu}}\setm \infi{P_{\alpha}}\subs \{s',t',b'\}$.

So we completed the proof of Lemma \ref{lm:alg}.
\end{proof}

After this preparation we can describe the winning strategy of Alice. 
We divide the game into
stages. 

\noindent {\bf 
Stage 1.} 

Let $F\subs {[{\omega}]}^{2}$ be fixed such that 
$\ddeg_F(v)=4$ for each $v\in {\omega}$.
Write $F=\{f_\ell:\ell<{\omega}\}$.
For $\ell<{\omega}$
let Alice ask the pair $f_\ell$ in the $\ell^{\rm th}$ turn.

Stage 1 terminates after ${\omega}$ turns. 
Then $\deg_{P_{\omega}}(n)=4$ for each $n\in {\omega}$.
So $$\BBB{\omega}\cap \TTT{\omega}=\empt$$ and clearly 
$\SSS{\omega}\subs \TTT{\omega}$. 

\medskip

\noindent {\bf 
Stage 2.} This stage is divided into substages. 
Before Substage $i$ the players played ${\eta}_i={\omega}+n_i$ turns 
for some $n_i<{\omega}$. Let $n_0=0$.

\smallskip

\noindent {\bf 
Substage i.}

If $\BBB{{\eta}_i}=\empt$ or $\SSS{{\eta}_i}=\empt$, then there is 
no scorpion graph $G$ which extends $G_{{\eta}_i}$, 
so the game has terminated. 
Since 
$\<{\omega},P_{{\eta}_i}\>$ is locally finite, i.e.
$P_{{\eta}_i}\notin \sid{1}$, Alice wins.  

So we can assume that $\BBB{{\eta}_i}\ne \empt$ and $\SSS{{\eta}_i}\ne\empt$.
Let 
\begin{displaymath}
k_i=\min(\BBB{{\eta}_i}\cup \SSS{{\eta}_i}),
\end{displaymath}
and let 
\begin{displaymath}
{C_i}=\left\{\begin{array}{ll}
{\BBB{{\eta}_i}}&\text{if $k_i\in \SSS{{\eta}_i}$,}\\
{\SSS{{\eta}_i}}&\text{if $k_i\in \BBB{{\eta}_i}$.}
\end{array}\right.
\end{displaymath}

Fix   an  enumeration  $\{e^i_j:j<M_i\}$  of 
$U_{{\eta}_i}\cap [\{k_i\},C_i]$ for some $M_i\le {\omega}$, and 
let Alice play the pair $e^i_j$ in the turn ${\eta}_i+j$.

If $k_i\in \BBB{{\eta}_i+j}\setm \BBB{{\eta}_i+j+1}$, or 
$k_i\in \SSS{{\eta}_i+j}\setm \SSS{{\eta}_i+j+1}$,
then Substage $i$ terminates, and we put 
${\eta}_{i+1}={\eta}_i+j+1$.

If the Substage $i$  has not terminated for any $j<M_i$, then 
we declare that the Substage $i$ was the last substage and 
Stage 2 terminates. 

\medskip
Assume that  Stage 2 has not terminated after  Substage $i$
for any $i<{\omega}$.  Then 
 we have 
\begin{displaymath}
      \min(\BBB{{\eta}_i}\cup \SSS{{\eta}_i})\notin  
      \BBB{{\eta}_{i+1}}\cup \SSS{{\eta}_{i+1}},
\end{displaymath}
for each $i<{\omega}$.
So taking ${\eta}=\sup\{{\eta}_i:i<{\omega}\}$
we have 
\begin{displaymath}
      \BBB{{\eta}}\cup \SSS{{\eta}}\subs 
      \bigcap_{i<{\omega}}(\BBB{{\eta}_{i+1}}\cup \SSS{{\eta}_{i+1}})
      =\empt.
\end{displaymath}
Thus, there is no scorpion graph $G$ which extends $G_{{\eta}}$, 
so the game terminates after at most ${\eta}$ steps. 
Since $\<{\omega},P_{{\eta}}\>$ is locally finite, i.e.
$P_{{\eta}}\notin \sid{1}$, Alice wins.

\noindent {\bf 
Stage 3.}

Assume that  Stage 2 terminated after Substage $i$, and
we are after turn ${\eta}$. Let us observe that ${\eta}={\eta}_i+{\omega}$ if $M_i={\omega}$,
and ${\eta}={\eta}_i+M_i+1$ if $M_i<{\omega}$.  Moreover, $$\infi{P_{\eta}}\subs \{k_i\}.$$

We should distinguish two cases.

\noindent{\bf Case 1.} $k_i\in \BBB{{\eta}_i}$.

Then $k_i\in \BBB{\eta}$ as well, and 
$[\{k_i\}, \SSS{{\eta}}]\subs P_{\eta}$.

Let
\begin{displaymath}
S=\{s\in S_{\eta}: \{s,k_i\}\in N_{\eta}\}.
\end{displaymath}

Assume  $G$ is a   scorpion graph which extends
$G_{\eta}$. Then $s_G\in \SSS{{\eta}}$, and so 
$\{s_G,k_i\}\notin E_G$, and so $\{s_G,k_i\}\in N_{\eta}$. 
Since $k_i\in \BBB{\eta}$ and so $\deg_{N_{\eta}}(k_i)\le 1$, it follows that 
$|S|\le 1$. So $S=\{s_G\}$.

Thus, if $|S|\ne 1$, then there is no scorpion graph $G$ 
which extends $G_{\eta}$, and
so the game terminates after at most ${\eta}$ turns. 
Since $\infi{P_{\eta}}\subs \{k_i\}$ and so  
$P_{{\eta}}\notin \sid{2}$, Alice wins.  

So we can assume that  $S=\{s'\}$, and  we know that 
\begin{displaymath}
\tag{A}\text{if $G \ge G_{\eta}$ is a scorpion graph, then $s_G=s'$.}
\end{displaymath}

Thus, we can apply Lemma \ref{lm:alg}(a):
Alice can 
      play in such a way that the game terminates in 
      ${\beta}\le {\eta}+{\omega}+{\omega}+{\omega}$
      turns and $|\infi{P_{\beta}}\setm \infi{P_{\eta}}|\le 3$.
Since $\infi{P_{\eta}}\subs \{k_i\}$, 
we have $|\infi{P_{\beta}}|\le 4$ and so 
$P_{{\beta}}\notin \sid{5}$, Alice wins.  

So we completed the investigation of Case 1.

\noindent{\bf Case 2.} $k_i\in \SSS{{\eta}_i}$.

Then $k_i\in \SSS{\eta}$ as well, and 
$[\{k_i\}, \BBB{{\eta}}]\subs P_{\eta}$.
Consider the set  
\begin{displaymath}
B=\{b\in B_{\eta}: \{k_i,b\}\in E_{\eta}\}.
\end{displaymath}

Since $k_i\in \SSS{\eta}$, we have  $\deg_{E_{\eta}}(k_i)\le 1$,
and so $|B|\le 1$.

\noindent {\bf Subcase 2.1} $B=\empt$.

First observe that 
\begin{displaymath}
\tag{A'}\text{if $G\ge G_{\eta}$  is a scorpion graph, then $s_G=k_i$}.
\end{displaymath}
Indeed,  in this case $b_G\in \BBB{\eta}$ and so 
$B=\empt$ implies $\{k_i,b_G\}\in N_{\eta}$, which yields $k_i=s_G$.

Thus, we can apply Lemma \ref{lm:alg}(a):
Alice can 
      play in such a way that the game terminates in 
      ${\beta}\le {\eta}+{\omega}+{\omega}+{\omega}$
      turns and $|\infi{P_{\beta}}\setm \infi{P_{\eta}}|\le 3$.
Since $\infi{P_{\eta}}\subs \{k_i\}$, 
$|\infi{G_{\beta}}|\le 4$ and so 
$G_{{\beta}}\notin \sid{5}$, Alice wins.

So we completed the investigation of Subcase 2.1.

\noindent {\bf Subcase 2.2} $B=\{b'\}$ for some $b'\in \BBB {\eta}$.

First observe that 
\begin{displaymath}
\tag{A''}\text{if $G\ge G_{\eta}$  is a scorpion graph, then $b_G=b'$}.
\end{displaymath}
Indeed, since $t_G\notin \BBB{\eta}$ because 
 $\BBB{\eta}\cap \TTT {\eta}\subs \BBB {\omega}\cap \TTT {\omega}=\empt$,   
$B\ne \empt$ implies $k_i\ne s_G$, and so $b_G\in B$. Since $|B|=1$,
we have $b'=b_G$.

Thus, we can apply Lemma \ref{lm:alg}(b):
Alice can 
      play in such a way that the game terminates in 
      ${\beta}\le {\eta}+{\omega}+{\omega}+{\omega}$
      turns and $|\infi{P_{\beta}}\setm \infi{P_{\eta}}|\le 3$.
Since $\infi{P_{\eta}}\subs \{k_i\}$, 
we have $|\infi{G_{\beta}}|\le 4$ and so 
$P_{{\beta}}\notin \sid{5}$, Alice wins.

So we completed the investigation of Subcase 2.2.

Since there are no more cases,  we proved Theorem \ref{tm:scorpion}.
\end{proof}

\begin{theorem}\label{tm:we-degree}(1)
      For each natural number $n\ge 1$ 
      the monotone graph property  $D_n$
       is  $\isf$-hard on the vertex set ${\omega}$.
      
       \noindent 
      (2)
      For each natural number $m\ge 2$ 
      the monotone graph property  $C_m$
      is  $\isf$-hard on the vertex set ${\omega}$.
\end{theorem}

\begin{proof} (1)
If $F\subs {[{\omega}]}^{2}$ and $j\in {\omega}$,
write 
$\ddeg_{F}(j)=|\{i\in {\omega}: \{i,j\}\in F\}|$
and 
$\dgl_{F}(j)=|\{i<j: \{i,j\}\in F\}|$.

Fix $n\in {\omega}$. Let Bob play using the following strategy:
in the ${\alpha}$th step, if Alice asks the undetermined pair $e_{\alpha}=\{i,j\}$ with 
$i<j<{\omega}$, then 
Bob says "yes" iff either
\begin{enumerate}[(1)]
\item[(a)] $i\le n$ and $\ddeg_{E_{\alpha}}(i)<n$, or $j\le n$ and $\ddeg_{E_{\alpha}}(j)<n$, 
\end{enumerate}
or 
\begin{enumerate}[(1)]
\item[(b)] $\ddeg_{E_{\alpha}}(j)+\dgl_{U_{\alpha}}(j)= n$. 
\end{enumerate}

Assume that the game terminates after ${\beta}$ steps.

\begin{lemma}\label{lm:dmaxgen}
For each ${\alpha}\le {\beta}$, we have       
$\ddeg_{\gmax{G_{\alpha}}}(v)\ge n$ for each $v\in {\omega}$.
\end{lemma}

\begin{proof}
Fix $j\ge n$. Then $\ddeg_{E_{0}}(j)+\dgl_{U_{0}}(j)=0+j\ge n$. 
      By transfinite induction on ${\alpha}$  we can see that      
       we have $\ddeg_{E_{\alpha}}(j)+\dgl_{U_{\alpha}}(j)\ge n$
      for each ${\alpha}\le {\beta}$.
      
  Moreover,   for each $i<n$ we have 
      $\ddeg_{E_{\alpha}}(i)\ge min(n, \ddeg_{P_{\alpha}}(i))$, so if 
      $\ddeg_{E_{\alpha}}(i)<n$, then $\ddeg_{U_{\alpha}}(i)$ is infinite.

      Thus, $\ddeg_{\gmax{G_{\alpha}}}(v)\ge n$ for each $v\in {\omega}$ and ${\alpha}\le {\beta}$.
      \end{proof}

\begin{lemma}\label{lm:isolatd} 
      If $\ddeg_{E_{\alpha}}(v)\ge 1$ for each $v\in {\omega}$ for some ${\alpha}\le {\beta}$,
      then $P_{G_{\alpha}}\in \isf$.
      \end{lemma}
      
      \begin{proof}[Proof of the Lemma]
      $E_{\alpha}$ contains just finitely many edges which were obtained by applying rule (a).
      Let $J\in {\omega}$ such that ${[J]}^{2}$ contains all of these edges.
      Let $\{i,j\}\in E_{\alpha}\setm {[J]}^{2}$. Then 
      there is ${\gamma}<{\alpha}$ such that $e_{\gamma}=\{i,j\}$ and rule (b)
      was applied, i.e. 
      $\ddeg_{E_{\gamma}}(j)+\dgl_{U_{\gamma}}(j)=n$, and so $\dgl_{U_{\alpha}}(j)\le n$.
       
      Let $A=\{j: J<j \land \{i,j\}\in E_{\alpha}\text{ for some } i<j\}$.

      Consider the following graph $K=\<A,F\>$: $\{i,j\}\in F$ iff $i<j$ and $\{i,j\}\in U_{\alpha}$.
      Since $\dgl_{F}(j)\le n$, the chromatic number of $K$ is at most $n+1$, 
      and so there is $B\in {[A]}^{{\omega}}$ such that ${[B]}^{2}\cap F=\empt$, i.e. 
      ${[B]}^{2}\subs P_{\alpha}$.
 So $P_{\alpha}\in \isf$.
      \end{proof}
         
By Lemma \ref{lm:dmaxgen} $\gmax{G_{\beta}}$ has property $D_n$. Since the game has terminated,
$\gmin{G_{\beta}}=\<{\omega},E_{\beta}\> $ also has property $D_n$. Since 
$n\ge 1$, we can apply Lemma \ref{lm:isolatd} for ${\alpha}={\beta}$
to obtain $P_{\beta}\in \isf$, which proves 
 (1).

\smallskip

 (2) Let Bob use the strategy he applied  again property $D_{m-1}$. 

Assume that the game terminates after ${\beta}$ moves. 

 Then  for each ${\alpha}\le {\beta}$  we have $\ddeg_{\gmax{G_{\alpha}}}(v)\ge m-1$ for each $v\in {\omega}$
 by Lemma \ref{lm:dmaxgen}.
Thus, $\cc v{\gmax{G_{\beta}}}\ge m$ for each $v\in {\omega}$ and so 
$\gmax{G_{\beta}}$ has property $C_m$.

Since the game has terminated, $\gmin {G_{\beta}}=\<V, E_{\beta}\>$
has property $C_m$ as well.
Then $\ddeg_{G_{\beta}}(v)\ge 1$ for each $v\in {\omega}$ because $m\ge 2$.
Thus, by Lemma \ref{lm:isolatd}, $P_{{\beta}}\in \isf$.  Thus, Bob wins.
\end{proof}

\begin{theorem}\label{tm:we-connected}
      The monotone graph property "connected" is $\isf$-hard on ${\omega}$.
      \end{theorem}
      
      \begin{proof}

       Let      Bob play according to the following strategy.
      
      If $e_{\alpha}=\{i,j\}$, then Bob says yes iff
      both $A=\CC i{G_{\alpha}}$, the connected component of $i$ in 
      $G_{\alpha}$,  and  $B=\CC j{G_{\alpha}}$, 
      the connected component of $j$ in 
      $G_{\alpha}$,  are finite, $A\ne B$ and 
      $[A,B]\subs P_{\alpha}\cup \{e_{\alpha}\}$.
      
      Assume that the game terminates after ${\beta}$ steps.

      Assume first that $G_{\beta}$ contains an infinite connected component, 
      and  let 
      \begin{displaymath}
      {\gamma}=\min\{{\gamma}'\le {\beta}: 
      \text{$G_{{\gamma}'}$ contains an infinite connected component}\},
      \end{displaymath}
     and consider an infinite component $A$ of $G_{\gamma}$.

     Fix  $\{i,j\}\in {[A]}^{2}$, and  let ${\alpha}$ be the minimal ordinal such that 
     $\CC i{G_{{\alpha}}}=\CC j{G_{{\alpha}}}$. Then 
     ${\alpha}={\delta}+1<{\gamma}$
     and $[\CC i{G_{{\delta}}},\CC j{G_{{{\delta}}}}]\subs P_{\alpha}$ because 
     Bob declared an edge between $\CC i{G_{{\delta}}}$ and $\CC j{G_{{{\delta}}}}$ in the 
     ${\delta}^{\rm th}$ turn. 
     Thus, $\{i,j\}\in P_{{\alpha}}$. Since $\{i,j\}$ was arbitrary, hence 
     ${[A]}^{2}\subs P_{{\gamma}}\subs P_{\beta}$. 
Thus, $P_{{\beta}}\in \isf$.

So we can assume that every connected component of $G_{\beta}$ is finite. 
So $G_{\beta}$ is not connected. Since the game has terminated, $\gmax{G_{\beta}}$
is not connected as well. 
So ${\omega}$ has a partition ${\omega}=X_0\cup X_1$ such that 
$(E_{\beta}\cup U_{\beta})\cap[X_0,X_1]=\empt.$
Pick $v_i\in X_i$ and let $K_i=\CC{v_i}{G_{\beta}}$ for $i<2$.
Since $K_i\subs X_i$, $[K_0,K_1]\subs N_{\beta}$.

But $[K_0,K_1]$ is finite, so let ${\alpha}$ be the maximal ordinal such that 
$e_{\alpha}\in [K_0,K_1]$. Let $\ell_i=e_{\alpha}\cap K_i$ for $i<2$. Then 
\begin{displaymath}
e_{\alpha}\in [\CC{\ell_0}{G_{\alpha}},\CC{\ell_1}{G_{\alpha}}]\subs [K_0,K_1]\subs
N_{\alpha}\cup \{e_{\alpha}\}, 
\end{displaymath}
and so
Bob declared that $e_{\alpha}$ is an edge. Contradiction, which proves that it is not possible that 
every connected component of $G_{\beta}$ is finite.

Thus, we proved $P_{\beta}\in \isf$.
\end{proof}

\begin{theorem}\label{tm:we-degree-2}
      For each natural number $n\ge 1$ 
      the monotone graph property 
      \begin{displaymath}
      \text{"{\em $G$  contains $K_{1,n}$}"  
      }
      \end{displaymath}
      is  $\mc F_n$-hard on the vertex set ${\omega}$,
      where 

      \begin{displaymath}
      \mc F_n=\{E\subs [{\omega}]^2:\exists B\in {[{\omega}]}^{n} \text{ such that } 
      {[{\omega}]}^{2}\setm  E\subs {[B]}^{2}. \}
      \end{displaymath}

\begin{proof}
Let Bob play according to the following greedy strategy:
in the ${\alpha}$th step if $e_{\alpha}=\{i,j\}$ then 
Bob says yes iff $\ddeg_{G_{\alpha}}(i)<n-1$ and $\ddeg_{G_{\alpha}}(j)<n-1$.

Assume that the game terminates after ${\beta}$ steps. 

Clearly $G_{\beta}$ does not contain $K_{1,n}$. 
Since the game has terminated, $\gmax{G_{\beta}}$ can not  contain $K_{1,n}$ as well. 

Let $A=\{v\in {\omega}:\ddeg_{E_{\beta}}(v)=n-1\}$
and $B=\{v\in {\omega}: \ddeg_{E_{\beta}}(v)\le n-2\}$.
Then $N_{\beta}\cap {[B]}^{2}=\empt$. 
Thus, $|B|\le n$ or $\ddeg_{E_{\beta}\cup U_{\beta}}(b)\ge |B|-1\ge n$ for each $b\in B$.

Moreover, 
$U_{\beta}\subs {[B]}^{2}$,
or $\{i,j\}\in U_{\beta}\cap {[A,{\omega}}]$ with 
$i\in A$
implies $d_{E_{\beta}\cup U_{\beta}}(i)\ge  n$.
Thus, $P_{\beta}\in \mc F_n$.

Thus, Bob wins.
\end{proof}

\end{theorem}

\bibliographystyle{plain}

\end{document}